
\input amstex
\input pictex
\documentstyle{amsppt}
\magnification\magstep1

\define\vo{\text{\rm vol}}

\topmatter
\title
Equivalence of Geometric and Combinatorial\\
Dehn Functions
\endtitle
\author
Jos\'e Burillo
\endauthor
\abstract
In this paper it is proved that if a finitely presented group
acts properly discontinuously, cocompactly and by isometries
on a simply connected Riemannian manifold, then the two
Dehn functions, of the group and of the manifold, respectively,
are equivalent.
\endabstract
\address
Dept. of Mathematics, University of Utah, Salt Lake City,
Utah 84112
\endaddress
\email
burillo\@math.utah.edu
\endemail
\endtopmatter

\document

\heading
1. Dehn functions and their equivalence
\endheading

Let $X$ be a simply connected 2-complex , and let $w$ be an
edge circuit in
$X^{(1)}$. If $D$ is a van Kampen diagram for $w$ (see \cite{5}), then
the area of $D$ is defined as the number of 2-cells on $D$, and the
area of $w$, $a(w)$, is defined as the minimum of the areas of all van
Kampen diagrams for $w$. Then the Dehn function of $X$ is defined as
$$
\delta_X(n)=\max a(w),
$$
where the maximum extends to all loops $w$ of length $l(w)\le n$.

Given two functions $f$, $g$ from $\Bbb{N}$ to $\Bbb{N}$ (or, more
generally, from $\Bbb{R}^+$ to $\Bbb{R}^+$), we say that $f\prec g$ if
there exist positive constants $A$, $B$, $C$, $D$, $E$ such that
$$
f(n)\le Ag(Bn+C)+Dn+E.
$$

Two such functions are called equivalent (denoted $f\equiv g$) if
$f\prec g$ and $g\prec f$. The importance of Dehn functions is given
by the fact that they are invariant by quasi-isometries: when one
considers the 1-skeleton of the complex as a metric space with the
path metric, where every edge has length one,
two complexes with quasi-isometric 1-skeleta have
equivalent Dehn functions (see \cite{1}).

Let $G$ be a finitely presented group, and let $\Cal{P}$ be a finite
presentation for $G$. Let $K=K(\Cal{P})$ be the 2-complex associated
with $\Cal{P}$, i.e. the 2-complex with a single vertex, an oriented
edge for every generator of $\Cal{P}$, and a 2-cell fer every relator,
attached to the edges according to the spelling of the relator. Then
the Dehn function of $\Cal{P}$ is, by definition, the Dehn function
$\delta_{\tilde K}$ of the universal covering of $K$. The fact that
two finite presentations $\Cal{P}$ and $\Cal{Q}$ for the same group
give 2-complexes $\widetilde{K}(\Cal{P})$ and
$\widetilde{K}(\Cal{Q})$ with quasi-isometric 1-skeleta,
and hence equivalent Dehn functions, leads
to the concept of Dehn function of the group $G$, meaning the
equivalence class of the Dehn function of any of its
presentations. An extensive treatment of Dehn functions of
finitely presented groups is given in \cite{4}.

A closely related definition can be formulated in the context of
riemannian manifolds, dating back to the isoperimetric problem for
$\Bbb{R}^n$ of Calculus on Variations. Given a Lipschitz loop $\gamma$
in a simply connected riemannian manifold $M$, we call the area of
$\gamma$ to the minimum of the areas of all Lipschitz discs bounding
$\gamma$. Then, clearly, we can define the Dehn function of $M$ as
$$
\delta_M(x)=\max_{l(\gamma)\le x}\text{area}(\gamma).
$$

It is natural to consider the question of whether the Dehn functions
of a simply connected riemannian manifold $M$ and of a finitely presented
group $G$ acting properly discontinuously and cocompactly on $M$
agree. The fact that they effectively agree has been implicitly
assumed in the literature, while no proof has been given. A closely
related statement is given in \cite{2,~{\rm Theorem 10.3.3}},
applying to this setting
the Deformation Theorem of Geometric Measure Theory (\cite{3,~4.2.9}
and \cite{7}),
and which provides the basis of the Pushing Lemma below. This paper is
devoted to provide a complete and detailed proof of the fact that the
two Dehn functions are equivalent. It is known to the author that
M. Bridson has an alternate, unpublished proof for the same result. The
author would like to thank Professor S. M. Gersten for his
encouragement and his useful remarks.

The statement of the theorem is as follows:

\proclaim{Theorem 1.1} Let $M$ be a simply connected riemannian
manifold, and let $G$ be a finitely presented group acting properly
discontinuously, cocompactly and by isometries on $M$. Let $\tau$ be a
$G$-invariant triangulation of $M$. Then the three following Dehn
functions are equivalent:
\roster
\item the Dehn function $\delta_G$ of any finite presentation of $G$,
\item the Dehn function $\delta_{\tau^{(2)}}$ of the 2-skeleton of
$\tau$, and
\item the Dehn function $\delta_M$ of $M$.
\endroster
\endproclaim

The fact that $\delta_G$ and $\delta_{\tau^{(2)}}$ are equivalent is
clear: since $G$ acts cocompactly on $\tau$, there is a
quasi-isometry between $\tau^{(1)}$ and the 1-skeleton of
$\widetilde{K}(\Cal{P})$ for
any presentation $\Cal{P}$ of $G$, and the equivalence follows from
the results in \cite{1}. We will concentrate on the proof of the
equivalence between $\delta_{\tau^{(2)}}$ and $\delta_M$, and the
arguments will be mainly geometric, trying to relate the lengths and
areas of loops and discs in $M$ with those included in the
triangulation $\tau$. The first step in this direction is given by the
Pushing Lemma, a complete analog of the Deformation Theorem in
Geometric Measure Theory and already stated and proved, in a slightly
different way, in \cite{2,~{\rm Theorem 10.3.3}}, whose proof we will
follow closely.

\heading
2. Technical Lemmas
\endheading

The main tool for the proof of the equivalence of the two Dehn
functions is the Pushing Lemma, which allows to relate an arbitrary
Lipschitz chain in $M$ with another chain which is included in the
corresponding skeleton of $\tau$. The fact that we must follow the
variation of volume of the chain and prevent its excessive growth is
what makes the argument more complicated, since projection from an
arbitrary point would lead to arbitrarily large growth of this
volume. Techniques from measure theory assure the existence of a
center of projection which is far enough from the chain, and hence
provides control on the growth of the volume.

\proclaim{Lemma 2.1 (Pushing Lemma)} Let $M$, $G$ and $\tau$ be as
above. Then there exist a constant $C$, depending only on $M$ and
$\tau$, with the following property:
Let $T$ be a Lipschitz $k$-chain in $M$, such that $\partial T$ is
included in $\tau^{(k-1)}$. Then there exist another Lipschitz
$k$-chain $R$, with $\partial R=\partial T$, which is included in
$\tau^{(k)}$, and a Lipschitz $(k+1)$-chain $S$, with $\partial
S=T-R$, and such that
$$
\vo_k(R)\le C\vo_k(T)\qquad \text{and}\qquad
\vo_{k+1}(S)\le C\vo_k(T).
$$
\endproclaim

The only difference with the statement given in \cite{2} is the
fact that their statement for cycles can be extended to chains, since
the boundary of the chain is not modified, being included in the
$(k-1)$-skeleton. A statement for
cycles is not sufficient, since this lemma will be applied to the
discs as well as to the loops, and the fact that $\partial T=\partial
R$ is crucial in the proof of the main theorem.

\demo{Proof}
The proof will proceed by descending induction on the skeleta of
$\tau$. So assume that $T$ is included in $\tau^{(i)}$ but not in
$\tau^{(i-1)}$, for $i>k$. We want to proceed simplex by simplex,
choosing an appropriate point in each simplex and projecting radially
from this point the chain $T$ to the boundary of the simplex. The claim we
will prove is: there exists a constant $C$ such that for every simplex
there exists a point such that projecting $T$ from this point to the
boundary of the simplex does not increase the volume of the chain by
more than a multiplicative factor $C$.

To simplify the computations, since the action of the group is
cocompact, we can change equivariantly the riemannian metric of $M$ to
assume that every simplex is the unit euclidean simplex of side
length one. Let $\sigma$ be an $i$-simplex of $\tau$,
let $O$ be the barycenter of $\sigma$, and
let $r$ be a positive number such that the ball of center $O$ and
radius $3r$ is included in the interior of $\sigma$. Let $B$ be the
ball of center $O$ and radius $r$, let $u\in B$, and let $B_u$ be the
ball of center $u$ and radius $2r$. Clearly $B\subset B_u$, for all
$u$. Let $\pi_u$ be the radial projection with center $u$ of
$B_u\setminus \{u\}$ onto $\partial B_u$. Let $Q=T\cap\sigma$.
We want to see that there exists a constant $v_0$, independent of $T$
and $\sigma$, and there exists $u\in B$, clearly dependent on $T$,
with
$$
\vo_k(\pi_uQ)\le v_0\,\vo_k(Q).
$$

\beginpicture
\setcoordinatesystem point at -181 110
\setplotarea x from -181 to 100, y from -85 to 100
\circulararc 360 degrees from 15 0 center at 0 0
\circulararc 240 degrees from -7 -37 center at -7 -7
\setlinear
\plot 0 100  87 -50  -87 -50  0 100 /
\setquadratic
\setdashes <3pt>
\plot -33 8  -15 -15  -7 -37 /
\setsolid
\plot -50 14  -40 13  -33 8 /
\plot -7 -37  -5 -43  0 -50 /
\circulararc 120 degrees from -33 8 center at -7 -7
\put {$\sigma$} [lB] at 18 75
\put {$O$} [lB] at 1 1
\put {$u$} [l] at -5 -7
\put {$B$} [lt] at 0 -17
\put {$B_u$} [lt] at 15 -30
\put {$Q$} [rt] at -16 -16
\put {$\pi_u Q$} [rt] at -30 -30
\put {.} at 0 0
\put {.} at -7 -7
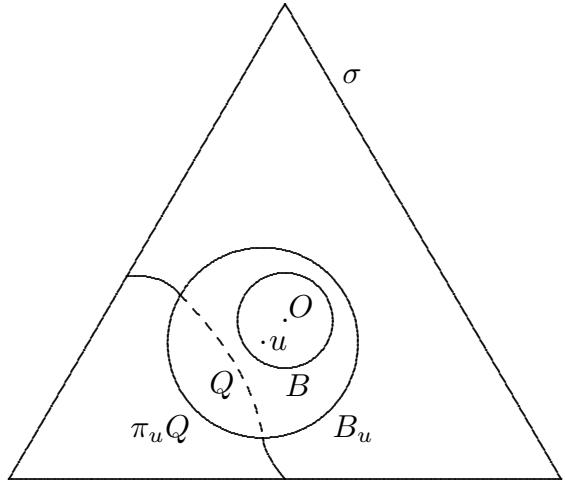
\put {Figure 1: \sl Projecting $Q$ to the boundary of $B_u$.} [B] at
0 -70
\endpicture

For every positive real number $v$ define
$$
A_v=\{u\in B\,|\,\vo_k(\pi_uQ)>v\,\vo_k(Q)\}
$$
and let $\alpha(v)=m_i(A_v)$, where $m_i$ is the $i$-dimensional
Lebesgue measure. We want to prove that
$$
\lim_{v\to\infty}\alpha(v)=0,
$$
so it will be enough to choose $v_0$ such that $\alpha(v_0)<m_i(B)$ to
have $B\setminus A_{v_0}$ nonempty.

We have
$$
\split
\vo_k(\pi_uQ)&\le\vo_k(\pi_u(Q\cap B_u))+\vo_k(Q)\\
&\le\int_{Q\cap B_u}\left(\frac{2r}{||x-u||}\right)^k\,dx+\vo_k(Q),
\endsplit
$$
where the first term accounts for the volume obtained after
projecting, and the second term takes care of the possibility of $Q$
and $B_u$ being disjoint. Assume now that $\vo_kQ$ is nonzero (if
$\vo_kQ=0$ then $\vo_k(\pi_uQ)=0$). Then we have:
$$
\split
\alpha(v)\,v\,\vo_k(Q)&=v\,\vo_k(Q)\int_{A_v}du=\int_{A_v}v\,
\vo_k(Q)\,du\\
&\le\int_{A_v}\vo_k(\pi_uQ)\,du\le\int_B\vo_k(\pi_uQ)\,du\\
&\le\int_B\left(\int_{Q\cap B_u}\left(\frac{2r}{||x-u||}\right)^k\,dx
+\vo_k(Q)\right)\,du\\
&=(2r)^k\int_{Q\cap B_u}\int_B ||u-x||^{-k}\,du\,dx+
\vo_i(B)\vo_k(Q)\\
&\le(2r)^k\int_{Q\cap B_u}dx\int_{B(O,3r)}||u||^{-k}\,du+
\vo_i(B)\vo_k(Q)\\
&\le K\vo_k(Q),
\endsplit
$$
where
$$
K=(2r)^k\int_{B(O,3r)}||u||^{-k}\,du+\vo_i(B).
$$
Observe that $K$ is finite and independent of $T$ and of $\sigma$.
The conclusion
is that $\alpha(v)v\le K$. Now, knowing $K$, we can find $v_0$
such that $K/v_0<m_i(B)$, and this $v_0$ is a constant
independent of $T$ and $\sigma$. We have found now $A_{v_0}$ with
strictly less measure than $B$, so we can pick a point in $B\setminus
A_{v_0}$ from which to project and make sure that the volume increases
at most by a multiplicative factor $v_0$.

The result of the above argument is the construction of another chain
$\pi_uQ$ which is far enough from $O$. We can now project radially
from $O$ to $\partial\sigma$, and the change of volume is bounded
since $\pi_uQ$ is at least at a distance $r$ from $O$. The combination
of this change of volume with $v_0$ gives the constant needed in this
precise skeleton. Combining all the constants from all the steps we
obtain the desired constant $C$.
Observe that these projections leave $\tau^{(i-1)}$ unchanged, so
clearly $\partial T$ is preserved by them.

The $(k+1)$-chain $R$ is obtained by joining every $x\in Q$ with
$\pi_ux$ with a segment. The volume of the piece of $R$ contained in
$\sigma$ is bounded then, as before, by
$$
(2r)^{k+1}\int_{Q\cap B_u}\frac{dx}{||x-u||^k},
$$
where the extra factor $2r$ is obtained from the direction of the
projection, since each segment has length bounded by $2r$. An argument
similar to the previous one shows that projecting from most points
in $B$ gives the right bound for the volume.  $\square$
\enddemo

The second lemma states that for a Lipschitz map, almost every point
in the target space has a finite number of preimages. It is a direct
consequence of the area formula for Lipschitz maps, and it will be
used for both loops and discs in the proof of Theorem 1.1.

\proclaim{Lemma 2.2} Let $T$ be a Lipschitz $k$-chain in $M$, where
$k\le \dim M$. Then the set of points in $M$ with infinite preimages
under $T$ has Hausdorff $k$-measure zero.
\endproclaim

\demo{Proof}
Let $\sigma_k$ be the standard closed $k$-simplex, and let
$$
E:\sigma_k\longrightarrow M
$$
be one of the simplices in $T$. Since $E$ is a Lipschitz map,
by Rademacher's Theorem (\cite{3,~3.1.6}) it is
differentiable almost everywhere (with respect to the Lebesgue
$k$-measure), so the Jacobian $J_k E(x)$ is well defined for almost
all $x\in \sigma_k$. For $y\in M$,
let $N(E,y)$ be the number of elements of
$E^{-1}(y)$, possibly infinite,
and denote by $m_k$ and $h_k$ the Lebesgue and Hausdorff
$k$-measures, respectively. Then the area formula for Lipschitz maps
(\cite{3,~3.2.3}) states that
$$
\int_{\sigma_k} |J_k E(x)|\,dm_k(x)=\int_M N(E,y)\,dh_k(y).
$$
Since $E$ is Lipschitz, we have that $|J_k E(x)|$ is bounded,
and since $\sigma_k$ has finite measure, the
integral on the left hand side is finite. So the set where $N(E,y)$ is
infinite cannot have positive Hausdorff $k$-measure, because the right
hand side would be infinite. $\square$
\enddemo

\heading
3. Proof of the first inequality
\endheading

We will now prove one of the two inequalities involved in the
proof of the equivalence of $\delta_M$ and
$\delta_{\tau^{(2)}}$. Namely,
$$
\delta_M\prec\delta_{\tau^{(2)}}.\tag3.1
$$
Let $\gamma$ be a Lipschitz loop in $M$, with length at most
$n$. Using the Pushing Lemma, we can construct a new loop $\eta$, of
length at most $Cn$, which is included in the 1-skeleton,
and the homotopy between $\gamma$ and $\eta$ has
area at most $Cn$.

The loop $\eta$ is not combinatorial, so we will construct a homotopy
between $\eta$ and a new loop $\zeta$ which will be combinatorial, and
this homotopy will be included in the 1-skeleton, so it will have area
zero. Assume that $\eta$ is parametrized by
$$
\eta:S^1\longrightarrow M.
$$
Since every edge $e$ in $\tau^{(1)}$ has positive Hausdorff
$1$-measure, by Lemma 2.2 we can choose a point $p_e$
in the interior of $e$ such that
the set $\eta^{-1}(p_e)$ is finite. Let
$$
\eta^{-1}(\{p_e\,|\,e\text{ edge of }\tau\})=\{\theta_1,\ldots,\theta_m\}
\subset S^1
$$
with
$$
0\le \theta_1<\ldots<\theta_m<2\pi.
$$
This gives us a partition of the circle $S^1$ into arcs
$[\theta_i,\theta_{i+1}]$, for $i=1,\ldots,m$ (where $\theta_{m+1}=\theta_1$),
such that for every $i$, one of the two following
situations must occur:
\roster
\item $\eta(\theta_i)=\eta(\theta_{i+1})$, or
\item $\eta(\theta_i)=p_{e_i}$ and $\eta(\theta_{i+1})=p_{e_{i+1}}$
where $e_i$ and $e_{i+1}$ are two edges with a vertex $v_i$ in common.
\endroster
In the first case, we can collapse $\eta([\theta_i,\theta_{i+1}])$ into
$\eta(\theta_i)$, and we can construct a new parametrization of $\eta$
where for every $i$, $\eta(\theta_i)$ and $\eta(\theta_{i+1})$ are different.
Observe that $\eta([\theta_i,\theta_{i+1}])$ is exactly the
concatenation of the two segments $[p_{e_i},v_i]$ and
$[v_i,p_{e_{i+1}}]$, although the map is not a homeomorphism.
Find a homotopy between $\eta$ and the loop $\zeta$ such that:
$$
\zeta(\theta_i)=\eta(\theta_i),
$$
and $\zeta$ maps $[\theta_{i},\theta_{i+1}]$ homeomorphically to the
concatenation of the segments $[p_{e_i},v_i]$ and $[v_i,p_{e_{i+1}}]$.
This homotopy is length-decreasing: call $l_i$ the length of the
concatenation of the two segments. Then:
$$
\int_{\theta_{i}}^{\theta_{i+1}}|\eta'(t)|\,dt\ge
\left|\int_{\theta_{i}}^{\theta_{i+1}}\eta'(t)\,dt\right|=l_i.
$$
So the length of $\zeta$ is at most the length of $\eta$, and hence at
most $Cn$. Choose $\zeta^{-1}(\tau^{(0)})$ now as the set of vertices
of a simplicial structure in $S^1$. There is only one situation now
which prevents the map $\zeta$ from being combinatorial: it is possible that
an edge in $S^1$ is mapped to a loop starting and ending in the same vertex,
but since this loop is contained in a single edge of $\tau$, it can be
contracted to the vertex. The result of this contraction is now
simplicial. Contracting the corresponding edges in $S^1$ we obtain
finally a combinatorial loop in $M$, of length at most $Cn$.

This combinatorial loop can be
filled combinatorially by at most $\delta_{\tau^{(2)}}(Cn)$
2-simplices in $\tau$. The conclusion is that
$$
\delta_M(n)\le\delta_{\tau^{(2)}}(Cn)+Cn,
$$
and $\delta_M\prec\delta_{\tau^{(2)}}$.

\heading
4. Construction of a simplicial disc
\endheading

To prove the reverse inequality to (3.1), we start with a combinatorial loop
$\gamma$ in the 1-skeleton of $\tau$, with length at most $n$. Let
$$
f:D^2\longrightarrow M
$$
be a disc in $M$ with boundary $\gamma$, and with area $a$. We want to
construct a van Kampen diagram for $\gamma$ and bound its area in terms
of $a$. The first step is, as before, to use the Pushing Lemma to find
a new disc (also denoted $f$) which is included in $\tau^{(2)}$, and
whose area is at most $Ca$.

Let $\sigma$ be an open 2-simplex of $\tau$. Again by Lemma 2.2,
we can choose a point $p\in\sigma$, such that $f^{-1}(p)$ is
finite. Let $X$ be a component of $f^{-1}(\sigma)$. If $X\cap
f^{-1}(p)=\varnothing$, then $f\big|_{\dsize X}$ can be pushed
into $\partial\sigma$ projecting radially from $p$, and the resulting
map has area zero. So assume that $X$ is a component of
$f^{-1}(\sigma)$ with $X\cap f^{-1}(p)\ne\varnothing$,
and note that there are only finitely many of these components.
We can also assume that $f\big|_{\dsize X}$ is surjective, since if it
is not we can again project radially from a point not in $f(X)$.
Observe that even though these two radial projections can destroy the
Lipschitz property for $f$, from now on we are only going to use the
Lipschitz condition on $f\big|_{\dsize X}$, for those $X$ on which $f$
has not been modified by any radial projection.

Our way to find a lower bound on the area of $f\big|_{\dsize X}$ will
be using the degree of $f\big|_{\dsize X}$.
Since $f\big|_{\dsize X}$ is differentiable almost everywhere,
we can define the degree of $f\big|_{\dsize X}$ at a
point $y\in f(X)$ as
$$
\text{deg}\,f\big|_{\dsize X} (y) = \sum_{x\in f^{-1}(y)}
\text{sign}\,J_2f(x)
$$
(see \cite{3,~4.1.26}). Moreover, the degree of $f\big|_{\dsize X}$
is almost constant in $f(X)$, so we can define the degree of
$f\big|_{\dsize X}$ as the value $d_X$ it achieves at almost every
$y\in f(X)$. The lower bound is given by the area formula for
Lipschitz maps: if $u$ is an integrable function respect to $m_2$, we
have (see \cite{3,~3.2.3}):
$$
\int_X u(x)|J_2f(x)|\,dm_2=\int_{\sigma} \sum_{x\in f^{-1}(y)\cap X}
u(x)\,dh_2,
$$
and taking $u(x)=\text{sign}\,Jf(x)$ we obtain:
$$
\gather
\text{area}\, f\big|_{\dsize X} =\int_X |J_2f(x)|\,dm_2\ge\left|\int_X
J_2f(x)\,dm_2\right|=\\
\left|\int_X \text{sign}\,J_2f(x)\,|J_2f(x)|\,dm_2\right|=
\left|\int_{\sigma}
\text{deg}\,f\big|_{\dsize X}\,dh_2\right|=\frac{\sqrt 3}2 |d_X|.
\endgather
$$

Our goal is to find a simplicial map
$$
g:D^2\longrightarrow \tau^{(2)}
$$
(with some simplicial structure in $D^2$) such that in $g\big|_{\dsize
X}$ only $|d_X|$ simplices are mapped by the identity to $\sigma$, and the
rest of $X$ is mapped to $\partial\sigma$. Then we will have that the
combinatorial area of $g$ is bounded the following way:
$$
\sum_X |d_X|\le \sum_X \frac2{\sqrt 3}\text{area}\left(f\big|_{\dsize
X}\right)\le\frac 2{\sqrt 3}Ca
$$
giving us the required bound. A technical result is needed since $g$
is not combinatorial, but only simplicial, and this result will be the
subject of the next section.

The first step in finding the map $g$ is to smooth the map
$f\big|_{\dsize X}$, to be
able to use differentiable techniques on it. Let $O$ be the barycenter
of $\sigma$, and choose $0<\epsilon<r$ such that:
$$
\varnothing\ne B(O,r-\epsilon)\subset B(O,r)\subset B(O,2r)
\subset B(O,2r+\epsilon)\subset\sigma,
$$
and let $U_1=f^{-1}(B(O,r))$ and $U_2=f^{-1}(B(O,2r))$. We have that
$\overline{U_1}\subset U_2\subset\overline{U_2}\subset X$. Choose
$\delta>0$ such that $B(x,\delta)\subset X$ for all $x\in U_2$, and
such that if $|x-y|<\delta$ then $|f(x)-f(y)|<\epsilon$, for all
$x,y\in X$. Let $\varphi$ be a $C^\infty$ bump function in $\Bbb{R}^2$
with support in $B(0,\delta)$, and with integral 1. Then,
for $x\in U_2$, we can construct the convolution
$$
f*\varphi(x)=\int_{B(x,\delta)}f(x-z)\varphi(z)\,dz,
$$
which is $C^\infty$ in $U_2$, and satisfies
$|f(x)-f*\varphi(x)|<\epsilon$ for all $x\in U_2$.
Also, if $f\big|_{\dsize X}$
was Lipschitz with constant $L$, then $f*\varphi$ is also
Lipschitz with the same constant: if $x,y\in U_2$,
$$
|f*\varphi(x)-f*\varphi(y)|\le|f(x-z)-f(y-z)|\int_{B(0,\delta)}
\varphi(z)\,dz\le L|x-y|.
$$
Choose now a Lipschitz function $\alpha$ on $X$ with values in
$[0,1]$ and equal to 1 in $U_1$ and to 0 outside $U_2$, and define
$$
\tilde f =\alpha (f*\varphi) + (1-\alpha)f\big|_{\dsize X}.
$$
Note that $\tilde f$ is defined only on $X$.
Then $\tilde f$ satisfies the following properties:
\roster
\item $|f(x)-\tilde f(x)|<\epsilon$ for all $x\in X$,
\item $\tilde f$ is smooth in $U_1$,
\item $\tilde f=f$ in $X\setminus U_2$,
\item $\tilde f$ is Lipschitz, and
\item $\text{deg}\,\tilde f=\text{deg}\,f\big|_{\dsize X}$.
\endroster
The first three properties are clear from the construction, and
property (4) holds because
$f\big|_{\dsize X}$, $f*\varphi$ and $\alpha$ are all
Lipschitz. To see that the degree is unchanged, since the degree is
almost constant, and $f\big|_{\dsize X}$
and $\tilde f$ agree outside $U_2$, we only
need to choose a point in $\sigma\setminus B(O,2r+\epsilon)$ for which
the degree is $d_X$ for both $f\big|_{\dsize X}$ and $\tilde f$.

We can now use Sard's Theorem (\cite{6})
to claim the existence of a regular
value for $\tilde f$ in $B(O,r-\epsilon)$ whose preimages are all in
$U_1$. Let $q$ be this regular value and let $p_1,\ldots,p_m$ be its
preimages. Let $V$ be an open disc with center $q$ such that
$\tilde f^{-1}(V)=V_1\cup\ldots\cup V_m$, where the $V_i$ are discs
around $p_i$, pairwise disjoint,
and such that $\tilde f\big|_{\dsize V_i}$ is a diffeomorphism. In
general, we will have that $m>|d_X|$, for which we will have to cancel
discs with opposite orientations. Assume $V_{m-1}$ and $V_m$ are
mapped to $V$ with opposite orientations. Choose $a\in\partial
V_{m-1}$ and $a'\in\partial V_m$ with $\tilde f(a)=\tilde f(a')$, and
join $a$ and $a'$ with a simple path $\lambda$ such that
$\tilde f(\lambda)$ is nullhomotopic in $\sigma\setminus V$,
which can be done because the map
$$
\tilde f:X\setminus \bigcup_{i=1}^{m}V_i\longrightarrow
\sigma\setminus V
$$
induces a surjective homomorphism of the fundamental groups.
Contracting $\tilde f(\lambda)$ we can assume $\tilde f(\lambda)$
is the constant path $\tilde f(a)$. Remove the discs
$V_{m-1}$ and $V_m$ and perform surgery along $\lambda$. The new
boundary thus created is mapped to $\partial V$ under $\tilde f$ by a
map from $S^1$ to itself of degree zero. Extend this map to a map from
$D^2$ to $S^1$ and attach it to $\tilde f$ along this boundary. For the
new map (which we will continue calling $\tilde f$),
the preimage of $q$ consists only of the points
$p_1,\ldots,p_{m-2}$. Repeating this process we will obtain a map
where now only the discs $V_1,\ldots,V_{|d_X|}$ are mapped to
$V$, and all with the same orientation.

Choose (temporarily) a sufficiently fine subdivision of $\tau$ such
that there is a 2-simplex $W$ in $V$, and let
$\rho_i={\tilde f}^{-1}(W)$.
Modify the map in $X$ by composing with the expansion of $W$ into all
$\sigma$.

\beginpicture
\setcoordinatesystem point at 0 78
\setplotarea x from 0 to 361, y from -103 to 70
\setdots <3pt>
\setlinear
\plot 40 68  160 68 /
\plot 30 51  170 51 /
\plot 20 34  180 34 /
\plot 10 17  190 17 /
\plot 0 0  180 0 /
\plot 10 -17  170 -17 /
\plot 20 -34  120 -34 /
\plot 30 -51  110 -51 /
\plot 40 -68  100 -68 /

\plot 0 0  40 68 /
\plot 10 -17  60 68 /
\plot 20 -34  80 68 /
\plot 30 -51  100 68 /
\plot 40 -68  120 68 /
\plot 60 -68  140 68 /
\plot 80 -68  160 68 /
\plot 100 -68  170 51 /
\plot 150 -17  180 34 /
\plot 170 -17  190 17 /

\plot 0 0  40 -68 /
\plot 10 17  60 -68 /
\plot 20 34  80 -68 /
\plot 30 51  100 -68 /
\plot 40 68  110 -51 /
\plot 60 68  120 -34 /
\plot 80 68  130 -17 /
\plot 100 68  150 -17 /
\plot 120 68  170 -17 /
\plot 140 68  180 0 /
\plot 160 68  190 17 /

\setsolid
\setquadratic
\plot 100 60  40 40  30 0  70 -60  100 -30  110 -10  130 0  160 10
      150 30  120 54  100 60 /

\setlinear
\plot 249 10  243 -5  255 -5  249 10 /
\plot 249 30  223 -15  275 -15  249 30 /
\plot 334 30  308 -15  360 -15  334 30 /

\arrow <7pt> [.25,.5] from 190 0 to 223 0
\arrow <7pt> [.25,.5] from 275 0 to 308 0

\plot 60 34  50 17  70 17  60 34 /
\plot 90 17  100 0  110 17  90 17 /
\plot 70 -17  60 -34  80 -34  70 -17 /

\put {$X$} [lB] at 139 44
\put {$\rho_1$} at 60 22.5
\put {$\rho_2$} at 100.5 11.5
\put {$\rho_3$} at 70 -28.5
\put {$W$} [lt] at 255 0
\put {$\sigma$} [lB] at 265 10
\put {$\sigma$} [lB] at 350 10
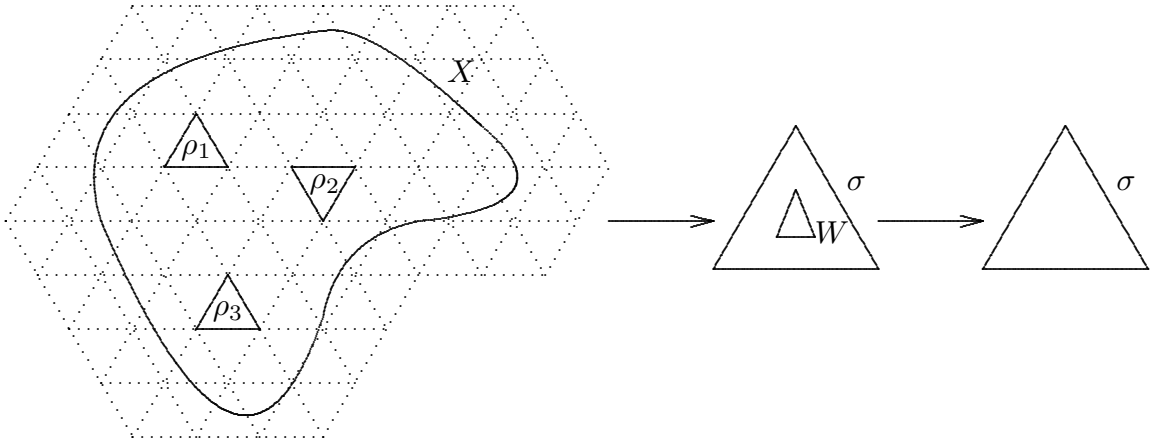
\put {Figure 2: \sl Making the map $f$ simplicial} [B] at 180.5 -88

\endpicture

After this process is done
for all $\sigma$, we obtain a map from $D^2$ to $\tau^{(2)}$, where
all the $\rho_i$ are sent homeomorphically to
2-simplices of $\tau$, and the rest is sent to the 1-skeleton of $\tau$.
To finish the construction
of $g$, find a simplicial structure on $D^2$ compatible with the
simplicial structure on the original loop $\gamma$ and which includes
all the $\rho_i$ obtained for all $\sigma$ as 2-simplices, and
approximate simplicially within $\tau^{(1)}$ the map $\tilde f$
relatively to all the $\rho_i$ and to $\gamma$.
The result is now simplicial, and the number of simplices sent by $g$
homeomorphically to 2-simplices in $\tau$ is
$$
\sum_X|d_X|\le\frac2{\sqrt 3}Ca.
$$
This map
is not a van Kampen diagram yet, since is only simplicial, and that is
the subject of the next section.

\heading
5. Degenerate Dehn functions
\endheading

Recall from the definition of van Kampen diagram that the map
is required to be combinatorial, i.e. every open cell is mapped
homeomorphically onto an open cell of the target.
It would be useful to extend this definition to
maps which are only simplicial, as the one obtained in the previous
section. This leads to the following definitions:

{\bf Definition 5.1:} Let $K$ be a simplicial 2-complex, and let
$w$ be a simplicial loop in $K^{(1)}$. A degenerate van Kampen diagram for
$w$ is a simplicial map from a planar contractible 2-complex $D$, with
some simplicial structure,
to $K$, such that the map restricted to the $\partial D$ is $w$. The length
of $w$ is defined as the number of 1-simplices on $S^1$ which are
mapped homeomorphically to 1-simplices of $K$, and similarly the area
of the degenerate van Kampen diagram is the number of 2-simplices of
$D$ which are mapped homeomorphically to 2-simplices of $K$. Then,
given a path $w$, its area is defined as the minimum of the areas of
all degenerate van Kampen diagrams for $w$. And the
degenerate Dehn function of $K$ is defined as:
$$
\delta_K^{\text{deg}}(n)=\max_{l(w)\le n} \text{area}(w).
$$

In the context of the previous section, we have proved the inequality
$$
\delta^{\text{deg}}_{\tau^{(2)}}\prec\delta_M,
$$
since the map $g$ constructed in section 4 is a degenerate van Kampen
diagram. The result that finishes the proof of the main theorem is the
following:

\proclaim{Theorem 5.2} Let $K$ be a simplicial 2-complex. Then,
$$
\delta_K\equiv\delta_K^{\text{deg}}.
$$
\endproclaim

\demo{Proof} One of the inequalities is obvious: let $w$ be a
simplicial loop in $K$. Modify the simplicial structure on $S^1$,
taking any 1-simplex which maps into a vertex and collapse it. This
produces a combinatorial loop, which can be filled with a
combinatorial disc, which is in particular simplicial. Then,
$$
\delta_K^{\text{deg}}\prec\delta_K.
$$
For the opposite inequality, let $w$ be now a combinatorial loop. We
can fill it with a simplicial map $f$ from the disc to $K$,
and the only thing we need to do
is construct a combinatorial disc with smaller area. Choose a vertex
of $K$, and let $L$ be a connected component of $f^{-1}(v)$. Change
the disc $D^2$ by
\roster
\item collapsing $L$ to a point, and
\item every 2-simplex with a face adjacent to $L$, but not in $L$,
has to be sent to an edge of $K$ adjacent to $v$. Collapse all these
simplices to edges.
\endroster

\proclaim{Lemma 5.3} The result of the collapsing indicated above is a
planar contractible simplicial complex with some 2-spheres attached to a
vertex.
\endproclaim

\demo{Proof} Attach a 2-cell $e$ to the boundary of $D^2$ to obtain a
finite cellular structure on $S^2$. Assume $L$ is contractible. Then
clearly the result of the collapsing will be a 2-manifold since every
edge is still adjacent to only two faces, and the star of every vertex
(including the vertex obtained in the collapsing) is an open disc. An
easy count of vertices, edges and faces gives the Euler characteristic
equal to 2, and hence the result is a 2-sphere. Eliminating the cell
$e$ we obtain a planar contractible simplicial complex.

If $L$ is not contractible, then the same argument can be applied to
every connected component of $S^2\setminus L$, obtaining a wedge of
2-spheres, and only one of them contains the cell $e$. Again
eliminating $e$ we obtain a planar contractible complex with some
spheres attached to a vertex. $\square$
\enddemo

Excising the spheres we obtain a
new degenerate van Kampen diagram, and since $w$ is a combinatorial
map, no edge of the boundary has been collapsed, and then the map on
the boundary is still $w$.
Doing this for all connected components of $f^{-1}(v)$, and for
every $v$, we obtain a new map from some planar contractible 2-complex
into $K$ which is now combinatorial, i.e. a (nondegenerate) van Kampen
diagram. Observe that this process
cannot increase the area, but can only decrease it when the 2-spheres
are cut off. Then we produced a van Kampen diagram for $w$ with
smaller area than the original degenerate van Kampen diagram. This
proves the inequality
$$
\delta_K\prec\delta_K^{\text{deg}}.
$$
$\square$
\enddemo

\enddocument